# Stably weakly shadowing transitive sets and dominated splittings

Dawei Yang

March 10, 2010

**Abstract**

We prove that for any $C^1$-stably weakly shadowing transitive set $\Lambda$, either $\Lambda$ is a sink or a source, or $\Lambda$ admits a dominated splitting.

## 1 Introduction

Shadowing properties have the physical meanings: even small errors occur at each iteration, one can track the resulting pseudo orbit by a true orbit for a stable (hyperbolic) system. A generalization of the classical shadowing called *weakly shadowing* was introduced by [5]. The weakly shadowing property is $C^0$ and $C^1$ generic by [5, 6]. One of the problems is to characterize the stably weakly shadowing diffeomorphisms. The following conjecture is given by S. Gan:

**Conjecture.** *A diffeomorphism $f$ is $C^1$-stably weakly shadowing if and only if $f$ is tame.*

Let's be more precise. Let $M$ be a compact $C^\infty$ Riemannian manifold without boundary. Let $\mathrm{Diff}^1(M)$ be the space of $C^1$ diffeomorphisms of $M$, for "space", we mean there is a usual $C^1$ metric defined on $\mathrm{Diff}^1(M)$. Given $f \in \mathrm{Diff}^1(M)$. For any $x \in M$, $\mathrm{Orb}_f(x) = \{f^n(x)\}_{n \in \mathbb{Z}}$ is the orbit of $x$ with respect to $f$. For $\varepsilon > 0$, $\{x_n\}_{n \in \mathbb{Z}}$ is called an $\varepsilon$-*pseudo orbit* if $d(f(x_n), x_{n+1}) < \varepsilon$ for any $n \in \mathbb{Z}$. One say that $f$ has the weakly shadowing property if for any $\varepsilon > 0$, there is $\delta > 0$ such that for any $\delta$-pseudo orbit $\{x_n\}_{n \in \mathbb{Z}}$, there is $x \in M$ such that $\{x_n\}_{n \in \mathbb{Z}} \subset B(\mathrm{Orb}(x), \varepsilon)$, where $\mathrm{Orb}(x)$ is the orbit of $x$. $f$ is called $C^1$-*stably weakly shadowing* if there is a neighborhood $\mathcal{U}$ of $f$ such that any $g \in \mathcal{U}$ is weakly shadowing. $f$ is *tame* if there is a neighborhood $\mathcal{U}$ of $f$ such that any $g \in \mathcal{U}$ has only finitely many chain recurrent classes[1]. Gan's conjecture is true if $\dim M = 2$ [13]. There is no answer on Gan's conjecture for higher dimensional diffeomorphisms.

If we focus on the "local" case, we can show that every $C^1$-stably weakly shadowing transitive set admits a dominated splitting. Let $\Lambda$ be a compact invariant set of $f$. We say that $\Lambda$ has the *weakly shadowing property* if for any $\varepsilon > 0$, there is $\delta > 0$ such that for any $\delta$-pseudo orbit $\{x_n\}_{n \in \mathbb{Z}} \subset \Lambda$, there is $x \in M$ such that $\{x_n\}_{n \in \mathbb{Z}} \subset B(\mathrm{Orb}(x), \varepsilon)$. For any set $U$, define $M_U(f) = \cap_{n \in \mathbb{Z}} f^n(U)$ to be the maximal invariant set of $f$ in $U$. For a compact invariant set $\Lambda$, we say that $\Lambda$ has $C^1$-*stably weakly shadowing property* if there are a compact neighborhood $U$ of $\Lambda$ and a $C^1$ neighborhood $\mathcal{U}$ of $f$ such that $M_U(g)$ has the weakly shadowing property for any $g \in \mathcal{U}$. For an invariant set $\Lambda$, if there are an invariant splitting of the tangent bundle $T_\Lambda M = E \oplus F$, together with two constants $C > 0$ and $\lambda \in (0, 1)$ such that $\|Df^n|_{E(x)}\| \|Df^{-n}|_{F(f^n(x))}\| \leq C\lambda^n$ for any $n \in \mathbb{N}$ and $x \in \Lambda$, one say that $\Lambda$ has a $((C, \lambda)$-*)dominated splitting*, $\dim E$ is called the *index* of this dominated splitting. The notion of dominated splitting is much weaker than the notion of hyperbolic splitting. But dominated splitting is also a robust property and it's an important mechanism for many dynamical phenomena. Recall that a compact invariant set $\Lambda$ is called *transitive* if $\omega(x) = \Lambda$ for some $x \in \Lambda$.

**Theorem A.** *For a transitive set $\Lambda$, if $\Lambda$ has the $C^1$-stably weakly shadowing property, and if $\Lambda$ is neither a hyperbolic sink nor a hyperbolic source, then $\Lambda$ admit a dominated splitting.*

One should notice that [7] proved a similar result for homoclinic classes. Each homoclinic class is transitive from hyperbolic theories.

---
[1]This definition is given by S. Gan. Other definitions of "tame" is given by C. Bonatti and F. Abdenur [1]

# 2 The reduction of the problem and periodic linear cocycles

For a periodic point $p$ of $f \in \text{Diff}^1(M)$, we list all the eigenvalues of $Df^{\pi(p)}$ as $\{\lambda_1, \lambda_2, \cdots, \lambda_d\}$ which verifies
$$|\lambda_1| \leq |\lambda_2| \leq \cdots \leq |\lambda_d|.$$
One say that $p$ is an *almost source* if $|\lambda_1| = 1$; $p$ is an *almost sink* if $|\lambda_d| = 1$. The following lemma is from [7, Lemma 3.2]:

**Lemma 2.1.** *If $\Lambda$ is a $C^1$-stably weakly shadowing set of $f$, then there is a neighborhood $\mathcal{U}$ of $f$ and a neighborhood $U$ of $\Lambda$ such that $M_U(g)$ contains neither almost sinks nor almost sources for any $g \in \mathcal{U}$.*

In the work, the metric between compact sets is the Hausdorff methic. We use $d_H$ to denote the distance of the Hausdorff distance. The limits of compact sets are under the Hausdorff distance. The following lemma is on the limit of uniformly dominated splitting. One can see a proof in [3, Lemma 1.4].

**Lemma 2.2.** *Give $C > 0$ and $\lambda \in (0,1)$. If there is a sequence of diffeomorphisms $\{f_n\}$ and a sequence of compact sets $\{\Lambda_n\}$ such that $\Lambda_n$ is a compact invariant set of $f_n$ and $\Lambda_n$ admits a $(C, \lambda)$-dominated splitting of index $i$ with respect to $f_n$, then if $\Lambda = \lim_{n \to \infty} \Lambda_n$ exists, then $\Lambda$ admits a $(C, \lambda)$-dominated splitting of index $i$ with respect to $f$.*

Let $p_n$ be a periodic point of $f_n$ with period $\pi(p_n)$. One says that $\{p_n\}$ is *uniformly contracting at the period* if there are $C > 0$, $\lambda \in (0,1)$ and $\iota \in \mathbb{N}$ such that for any $p_n$ with $\pi(p_n) > \iota$, then
$$\prod_{i=0}^{[\pi(p_n)/\iota]-1} \|Df_n^\iota(f_n^{i\iota}(p_n))\| \leq C\lambda^{[\pi(p_n)/\iota]}.$$
One says that $\{p_n\}$ is *uniformly expanding at the period* if it is uniformly contracting at the period for $\{f_n^{-1}\}$.

One can extract the following lemma from [11] and [10, Lemma II.4, Lemma II.5].

**Lemma 2.3.** *Assume $\lim_{n\to\infty} f_n = f$. If $p_n$ is a periodic sink of $\{f_n\}$. Then either $\{p_n\}$ is uniformly contracting at the period, or there is a sequence of diffeomorphisms $\{g_n\}$ such that $\text{Orb}_{f_n}(p_n)$ is an almost sink of $g_n$ and $\lim_{n\to\infty} g_n = f$.*

*Similarly, one can get the similar results for sources.*

**Lemma 2.4.** *Let $\Lambda$ be a transitive set of $f$. If there is a sequence of diffeomorphisms $\{f_n\}$ such that*

- *$f_n$ has a periodic point $p_n$ and $\{p_n\}$ is uniformly contracting at the period.*
- *$\lim_{n\to\infty} \text{Orb}_{f_n}(p_n)$ exists, and it's a subset of $\Lambda$.*

*Then $\Lambda$ is a sink.*
*One can get the similar results for $f^{-1}$ similarly.*

*Proof.* One can see the proof of this lemma in [8, Lemma 3.2-Lemma 3.4]. $\square$

**Proposition 2.5.** *Assume that $\Lambda$ is a $C^1$-stably weakly shadowing set. Also assume that $\Lambda$ is transitive and $\Lambda$ is not a periodic orbit. Then there is a neighborhood $\mathcal{U}$ of $f$ and a neighborhood $U$ of $\Lambda$ such that any $g \in \mathcal{U}$ has neither sinks nor sources in $U$.*

*Proof.* We will prove this proposition by absurd. If the conclusion is not true, without loss of generality one can assume that there is a sequence of diffeomorphisms $\{f_n\}$ such that each $f_n$ has a periodic sink $p_n$ such that $\lim_{n\to\infty} \text{Orb}_{f_n}(p_n)$ exists and it is a subset of $\Lambda$. By Lemma 2.3, either $\{p_n\}$ is uniformly contracting at the period or there is a sequence of diffeomorphisms $\{g_n\}$ such that $\text{Orb}_{f_n}(p_n)$ is an almost sink of $g_n$ and $\lim_{n\to\infty} g_n = f$. Since $\Lambda$ is a $C^1$-stably weakly shadowing set, by Lemma 2.1, $\{p_n\}$ is uniformly contracting at the period. Then by Lemma 2.4, $\Lambda$ is a periodic sink. This contradicts to the fact that $\Lambda$ is not a periodic orbit. $\square$

Pugh [12] proved the following lemma which guarantee the existence of periodic orbits by $C^1$ perturbations:

**Lemma 2.6.** *For any $C^1$ neighborhood $\mathcal{U}$ of $f$ and for any non-periodic point $x$ with the property $x \in \omega(x)$, there are $N \in \mathbb{N}$ and $\varepsilon_0 > 0$ such that for any $\varepsilon \in (0, \varepsilon_0)$ and $n \geq N$, there is $g \in \mathcal{U}$ with the following properties:*

- $g(z) = f(z)$ for any $z \in M \setminus (\cup_{i=0}^{n} B(f^i(x), \varepsilon))$.

- *There is a periodic point $y$ of $g$ such that $y \in B(x, \varepsilon)$ and $\text{Orb}_g(y) \cap (M \setminus (\cup_{i=0}^{n} B(f^i(x), \varepsilon))) \subset \text{Orb}_f(x)$.*

We will prove the following folklore result:

**Lemma 2.7.** *Let $\Lambda$ be a transitive set of $f$. For any $\varepsilon > 0$ and for any neighborhood $\mathcal{U}$ of $f$, there is a periodic point $p$ of $g \in \mathcal{U}$ such that $d_H(\text{Orb}_g(p), \Lambda) < \varepsilon$.*

*Proof.* Since $\Lambda$ is transitive, there is $x \in \Lambda$ such that $\omega(x) = \Lambda$. For $\varepsilon > 0$, there is $N_1 \in \mathbb{N}$ such that $d_H(\{x, f(x), \cdots, f^n(x)\}, \Lambda) < \varepsilon/2$ for any $n > N_1$. For any neighborhood $\mathcal{U}$ of $f$, one can get two constants $N$ and $\varepsilon_0$ from Pugh's closing lemma (Lemma 2.6). Without loss of generality, one can assume that $N > N_1$ and $\varepsilon_0 > \varepsilon$. For $n > N$, by Pugh's closing lemma, there are $g \in \mathcal{U}$ and $y \in B(x, \varepsilon/2)$ such that

- $y$ is a periodic point of $g$;

- $g^i(y) \in B(f^i(x), \varepsilon/2)$ and $g^i(y) \in \text{Orb}_f(x)$ when $g^i(y) \in (M \setminus (\cup_{i=0}^{n} B(f^i(x), \varepsilon)))$.

From above properties one can check the conclusion of this lemma is true. □

As a corollary,

**Corollary 2.7.1.** *Let $\Lambda$ be a transitive set of $f$. There are a sequence of diffeomorphisms $\{f_n\}$ and a sequence of points $\{p_n\}$ such that $p_n$ is a periodic point of $f_n$ and $\lim f_n = f$ and $\lim \text{Orb}(p_n) = \Lambda$.*

Now we give the proof of the main theorem. Let $\Lambda$ be a $C^1$-stably weakly shadowing transitive set of $f$. Moreover, $\Lambda$ is neither a sink nor a source. We will prove the main theorem by absurd, i.e., we assume that $\Lambda$ doesn't admit any dominated splitting.

**Claim.** *Under above assumptions, $\Lambda$ is not a periodic orbit.*

*Proof.* If $\Lambda$ is a periodic orbit, then $\Lambda$ is neither an almost sink nor an almost source. So $\lambda_d > 1$ or $\lambda_1 < 1$. In the first case, since $\Lambda$ is not a source, one can get the dominated splitting on $\Lambda$; in the second case one can do similarly. □

Since $\Lambda$ is transitive, by Corollary 2.7.1, there is a sequence of diffeomorphisms $\{f_n\}$ and a sequence of points $\{p_n\}$ such that $p_n$ is a periodic point of $f_n$ and $\lim f_n = f$ and $\lim \text{Orb}(p_n) = \Lambda$. Since $\Lambda$ is not a periodic orbit, one has $\pi(p_n) \to \infty$ as $n \to \infty$.

Let $\Sigma = \coprod_{n \in \mathbb{N}} \{p_n, f(p_n), \cdots, f^{\pi(p_n)-1}(p_n)\}$. One can define a natural $d$-dimensional vector bundle $E$ on $\Sigma$ as following: for any $x \in \Sigma$, the fibre on $x$ is $T_x M$. For any $i \in [0, \pi(p_n) - 1] \cap \mathbb{N}$, we define $h(f_n^i(p_n)) = f_n^{i+1}(p_n)$ and $A|_{E(f_n^i(p_n))} = Df_n(f_n^i(p_n))$. Thus $\mathcal{A} = (\Sigma, h, E, A)$ is a *bounded large periodic system* as in [4, Section 2.1-Section 2.3]. Then by [4, Theorem 2.2],

- either there is an infinite subset $\Sigma' \subset \Sigma$ which is invariant by $h$ such that the periodic linear cocycle $\mathcal{A}' = (\Sigma', h, E|_{\Sigma'}, A)$ admits a dominated splitting, which means that there are an invariant splitting $E|_{\Sigma'} = E^{cs} \oplus E^{cu}$ with respect to $A$, together with two constants $C > 0$ and $\lambda \in (0,1)$ such that for any $x \in \Sigma'$, for any $n \in \mathbb{N}$, one has $\|A^n|_{E^{cs}(x)}\| \|A^{-n}|_{E^{cu}(h^n(x))}\| \leq C\lambda^n$.

- Or, there is a perturbation $\mathcal{B}$ of $\mathcal{A}$ and an infinite invariant subsets $\Sigma'$ of $\Sigma$ such that for any $x \in \Sigma$, all eigenvalues of $B(h^{\pi(x)-1}(x)) \circ B(h^{\pi(x)-1}(x)) \circ \cdots \circ B(x)$ are real, with same modulus.

Under the help of Franks' Lemma [9, 10], we can translate the above statement for diffeomorphims: by taking a subsequence if necessary, either there are constants $C > 0$ and $\lambda \in (0,1)$ such that there is a $(C, \lambda)$-dominated splitting on the orbit $\{p_n\}$; or there is a sequence of diffeomorphisms $\{g_n\}$ such that $\lim_{n\to\infty} g_n = f$ and $\mathrm{Orb}_{f_n}(p_n)$ is also a periodic orbit of $g_n$ and all eigenvalues of $Dg_n^{\pi(p_n)}(p_n)$ are all real, and with the same modulus.

For the first case, we can get that $\Lambda$ has a dominated splitting by Lemma 2.2. In the second case, one can get that either $\Lambda$ is a sink or a source; or by an arbitrarily small perturbation there is an almost sink or an almost source in an arbitrarily small neighborhood of $\Lambda$, which contradicts to the fact that $\Lambda$ is $C^1$-stably weakly shadowing.

**Acknowledgements.** *The author would like to thank the support of Jilin University.*

Dawei Yang, School of Mathematics, Jilin University, Changchun, 130012, P.R. China
Email: yangdw1981@gmail.com